\newif\ifjournal
\newif\ifpgfgen
\newif\ifanon

\journalfalse 
\pgfgenfalse 
\anonfalse   

\ifjournal
\documentclass[a4paper]{article}
\else
\documentclass[a4paper]{article}
\fi

\usepackage[l2tabu, orthodox]{nag}
\usepackage{amsmath,amsfonts,amssymb,amsthm}
\usepackage{booktabs}
\usepackage[group-separator={,}]{siunitx}
\usepackage{fix-cm}

\ifjournal
\usepackage[disable]{todonotes}
\else
\usepackage{todonotes}
\fi

\usepackage{subcaption}
\usepackage{natbib}
\usepackage{fullpage}
\usepackage{microtype}
\usepackage{hyperref}
\usepackage[noabbrev]{cleveref}

\providecommand{\keywords}[1]{\textit{Keywords}: #1}

\DeclareMathOperator*{\argmax}{arg\,max}

\ifpgfgen
\usepackage{tikz}
\usepackage{pgfplots}
\pgfplotsset{compat=1.13}
\usetikzlibrary{external}
\tikzsetexternalprefix{figures/}
\makeatletter
\renewcommand{\todo}[2][]{\tikzexternaldisable\@todo[#1]{#2}\tikzexternalenable}
\makeatother
\usepgfplotslibrary{statistics}
\else 
\usepackage{graphicx}
\usepackage{tikzexternal}
\fi

\theoremstyle{definition}
\newtheorem{mydef}{Definition}[section]
\newtheorem{myex}{Example}[section]
\newtheorem{myremark}{Remark}
\newtheorem{myresult}{Result}

\usepackage{tikzscale}
\graphicspath{{./img/}{./figures/}{.}}

\makeatletter
\def\input@path{{./img/},{./}}
\makeatother

\begin{document}
\title{Dynamic pricing in retail with diffusion process demand}
\ifanon
\author{Anonymized}
\else
\author{A.~N.~Riseth\thanks{\texttt{anriseth@gmail.com} ---
    Mathematical Institute, University of Oxford,OX2 6GG.}  }
\fi

\date{\today}

\maketitle

\begin{abstract}
  When randomness in demand affects the sales of a product,
  retailers use dynamic pricing strategies to maximize their
  profits.
  In this article, we formulate the pricing problem as a
  continuous-time stochastic optimal control problem, and find
  the optimal policy by solving the associated Hamilton-Jacobi-Bellman (HJB) equation.
  We propose a new approach to modelling the randomness in the
  dynamics of sales based on diffusion processes. The model assumes a
  continuum approximation to the stock levels of the retailer which
  should scale much better to large-inventory problems than the
  existing Poisson process models in the revenue management literature.
  The diffusion process approach also enables modelling of the demand
  volatility, whereas Poisson process models do not.

  We present closed-form solutions to the HJB equation when there is no
  randomness in the system. It turns out that the deterministic
  pricing policy is near-optimal for systems with demand uncertainty.
  Numerical errors in calculating the optimal pricing policy may, in
  fact, result in a lower profit on average than with the heuristic
  pricing policy.
\end{abstract}

\section{Introduction}\label{sec:intro}
Consider a monopolist retailer who wants to design a dynamic pricing policy
for a product over a given period, in order to maximize their total
revenue and minimize the cost associated with handling unsold items
at a given terminal time.
Two important components in the decision process are the ability to
take into account the uncertainty associated with future cost and demand and to
optimally adjust for new knowledge as it arrives.
We will illustrate how to address both components in this article,
focusing on large-inventory limits and multiplicative demand
uncertainty.
Assume the retailer sells large volumes of its product at
high frequency compared to the total pricing period. In such a setting, it is
appropriate to model the sales process in a continuum limit, both for
product volume and for time, similar to~\cite{kalish1983monopolist}.

In the revenue management literature, most efforts to model demand uncertainty in
continuous time have focused on Poisson
processes. See, for example, the overviews by~\cite{bitran2003overview}
or~\cite{aviv2012dynamic}.
We stress that in other communities, such as financial markets,
both diffusion processes and jump-diffusions are common for modelling
demand and spot-prices~\citep{benth2014stochastic}. The survey by~\cite{carmona2014survey}
gives an indication of how flexible more advanced models used for commodity
markets can be. With this article, we wish to inspire the revenue
management community to take advantage of this research when modelling
uncertainty in their own domain. For the remainder of this section,
however, we mainly focus on modelling in the revenue management literature.

A pure Poisson process assumption is not compatible with taking a
time-continuum limit for sales volume with demand uncertainty and may
be better suited for demand
modelling in industries with lower product sales volumes, such as
the airline and hotel industries.
\cite{maglaras2006dynamic},~\cite{schlosser2015dynamic1},
and~\cite{schlosser2015dynamic2}
propose pricing heuristics similar to this
article, by considering a deterministic model based on the asymptotic continuum
limit of Poisson processes.
To the author's knowledge, however,
few attempts have been made to unify demand uncertainty with the continuum
limit. For example,~\cite{raman1995optimal} and~\cite{wu2016dynamic} model
demand uncertainty as increments of a Brownian motion. As we will show
in this article, their approaches lead to demand processes that admit
negative sales, with a probability approaching $1/2$ over
infinitesimally small time periods. We believe this is an important
factor for why so little research has been done in this area.
This article proposes a different
approach to modelling demand uncertainty in order to remedy this.
In our approach, the parameters of the system are described
as diffusion processes that are solutions to
stochastic differential equations (SDEs).
This enables modelling of the demand volatility, which Poisson
processes do not.
Our proposed approach can be combined with parameter estimation
methods already used by retailers
and also extends naturally to multiple products.
Modelling the demand over time as a diffusion process has previously
been done by~\cite{chambers1992estimation} at macroeconomic scale
with UK national data. His focus was on the data assimilation aspect and was not
applied to a setting of optimal control.

The retailer's dynamic pricing policies are stochastic processes that
control the SDE which describes the depletion of stock. In this
article, we seek an optimal pricing policy that maximizes the expected value of
the profit over a given pricing period.
One way to find an optimal pricing policy is to solve an associated nonlinear partial
differential equation (PDE), known as the Hamilton-Jacobi-Bellman (HJB) equation
\citep{pham2009continuous}.
We provide closed-form solutions for the HJB equation in the
deterministic case, for both linear
and exponential demand functions. The solution identifies two pricing
regimes, one where the retailer maximizes their profits without
depleting the inventory and another
where the retailer aims to maximize the price and still deplete
its inventory.
\citet{xu2006monopolistic} have considered a continuous-time pricing
problem where the uncertainty is modelled by a geometric Brownian
motion. Their expected demand function is unbounded, with the result that the
optimal pricing strategy ensures that all stock is sold by the
terminal time. The demand functions we consider are bounded, and by
introducing a penalty on unsold stock we capture the pricing regime
change that does not arise in~\cite{xu2006monopolistic}.

In financial markets traders face a similar problem to
that presented in this article, known as the optimal execution, or liquidation,
problem. There, a trader tries to sell, or purchase,
a particular amount of an asset by a predetermined time. See, for example,
\citet{cartea2015algorithmic} for an overview of this problem.
Instead of controlling the price, the focus of the retailer pricing
problem, the trader directly controls how much
of the product to sell at a particular time. If the expected demand
model is invertible, the retailer's pricing problem can be
reformulated to control the expected amount of stock to sell at each
time. This reformulation is sometimes chosen in the revenue management
community as well, see, for example,~\cite{bitran2003overview}.
In this article we focus on the formulation that writes the expected
demand model in terms of the price.

By investigating the terms in the HJB equation,
we identify the cases where the deterministic-case solution is appropriate.
Potentially significant changes to the
pricing policy for the stochastic system are at the interface between
the two pricing regimes: far away from this interface the
deterministic pricing policy is near-optimal. For example, the
expected price path is decreasing when one takes into account
uncertainty, while it is not for the deterministic heuristic.
For a risk-neutral decision maker, however, the differences in profit
are insignificant for most cases that may be relevant in industry.

The article is structured as follows: In \Cref{sec:modelling}, we
describe the modelling of the system and compare the new parameter uncertainty
approach to the existing Brownian increments approach.
Then, a formulation of the pricing problem and the associated HJB
equation is given in \Cref{sec:decision_formulation}. We also propose
a method to estimate the multiplicative factor in our model, in
order to implement the pricing policy in practice.
The optimal pricing policy in the deterministic limit is covered in
\Cref{sec:deterministic_hjb}, and the comparison to the stochastic
system is shown in \Cref{sec:stochastic_hjb}.
Extensions to the problem, such as other models for uncertainty, and
risk aversion, are discussed in \Cref{sec:extensions}.
Finally, we conclude and
suggest avenues for further research in \Cref{sec:conclusion}.

\section{Modelling demand and uncertainty}\label{sec:modelling}
For a given, positive amount of initial stock of a product, we are interested in
modelling the product sales over some finite time period.
Assume that the initial quantity of stock is large, and that there is a
substantial volume sold over time periods that are small compared to the total
period of interest. These assumptions can apply to many products sold
by large retailers. For example, in the monopoly setting, this leads to
a continuum model similar to that of~\citet{kalish1983monopolist}.
For a given product, denote the amount of stock left at time $t$ by $S(t)$,
and let $q(a)$ represent product demand at price $a$, per unit
time.
For simplicity, we assume $q$ does not explicitly depend on time.
In the continuum limit, the change in  stock at time $t$
is thus
\begin{equation}
  d\hat{S}(\hat{t})=
  \begin{cases}
    -\hat{q}(\hat{a})\,dt & \textnormal{if } \hat{S}(\hat{t})>0,\\
    0& \textnormal{if }  \hat{S}(\hat{t})\leq 0.
  \end{cases}
\end{equation}
For a pricing policy
$\alpha(t)$, the remaining stock at time $t$ is then
\begin{equation}
  \hat{S}(\hat{t})=\hat{S}(0)-\int_0^{\hat{t}}\hat{q}(\hat{\alpha}(\hat{u}))\,d\hat{u}.
\end{equation}
In the remainder of the article,
we emphasise the dependence of remaining stock on a particular
pricing policy $\alpha(t)$ using the superscript $S^\alpha$.

\begin{figure}[htbp]
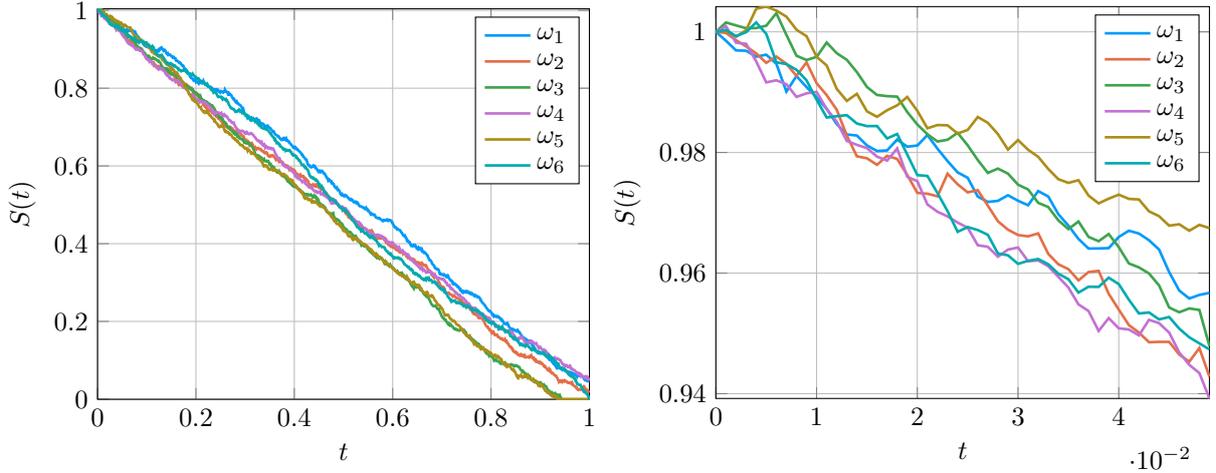

  \centering
  \includegraphics[width=0.5\textwidth,axisratio=1.25]{sde_brownian_sols}%
  \includegraphics[width=0.5\textwidth,axisratio=1.25]{sde_brownian_sols_zoom}
  \caption{Evolutions of stock $S^\alpha(t)$ for a constant
    demand forecast $q(a)=1$ with different sample paths $\omega_k$
    of a Brownian motion.
    The right hand figure shows that over small timescales, sales are
    often negative.
    All values are described in dimensionless quantities for
    simplicity, see Subsection~\ref{subsec:nondimensionalization}.
  }\label{fig:brownian_paths}
\end{figure}

At the start of the prediction period, it is not known exactly what
the demand will be at future times. We will now discuss how to
represent this uncertainty in the model. First, we note that the Brownian noise
approach in~\cite{raman1995optimal} and~\cite{wu2016dynamic} leads to negative
sales with probability tending to $0.5$ as the time period goes to zero.
Then, we propose a method that guarantees non-negative sales over all
time periods.
Let $W(t)$ denote a Brownian motion, and let $\sigma(t,s,a)$ be the volatility
in demand as a function of time, stock, and price. Then
one may say that the uncertainty in future
sales is due to the changes in $W(t)$, in the following sense.
\begin{equation}\label{eq:brownian_noise}
  d\hat{S}^{\hat{\alpha}}(\hat{t}) =
  \begin{cases}
    -\hat{q}(\hat{\alpha}(\hat{t}))\,d\hat{t} +
    \hat{\sigma}(\hat{t},\hat{S}(\hat{t}),\hat{\alpha}(\hat{t}))\,d\hat{W}(\hat{t})
    & \textnormal{if } \hat{S}(\hat{t})>0,\\
    0 &  \textnormal{if } \hat{S}(\hat{t})\leq0.
  \end{cases}
\end{equation}
In \Cref{fig:brownian_paths} we can see six realizations of
$S^\alpha(t)$ under
this model, using $q(a)=1$ and $\sigma(t,s,a)=0.05$. As we zoom in on
the sales paths, it is obvious that the stock often
increases over short time periods, corresponding to negative sales.

We now show that for any model with $\sigma>0$, the Brownian noise
model~\eqref{eq:brownian_noise} causes negative sales with high
probability.
For a given time period $[0,\Delta t)$, assume the demand and volatility are constant,
$q(a)=\tilde q>0$ and $\sigma(t,s,a)=\tilde \sigma>0$.
Then $S(\Delta t)=S(0)-\tilde q\Delta t + \tilde \sigma W(\Delta t)$.
Sales are negative in this period if $S(\Delta t)>S(0)$, that is when
$\tilde q\Delta t+\tilde \sigma W(\Delta t) < 0$.
The distribution of Brownian motion is normally distributed as
$W(\Delta t)\sim \mathcal N(0,\Delta t)$. So $W(\Delta t)$ is equal in
distribution to $\sqrt{\Delta t}Z$, where $Z\sim \mathcal N(0,1)$, and thus
\begin{align}
  \mathbb P(S(\Delta t)>S(0))
  &=\mathbb P(Z\leq -\sqrt{\Delta t}\,\tilde q / \tilde \sigma)
    \to 0.5^-,& \text{ as } \Delta t \to 0^+.
\end{align}

\begin{figure}[htbp]
  \centering
  \includegraphics[width=0.5\textwidth,axisratio=1.25]{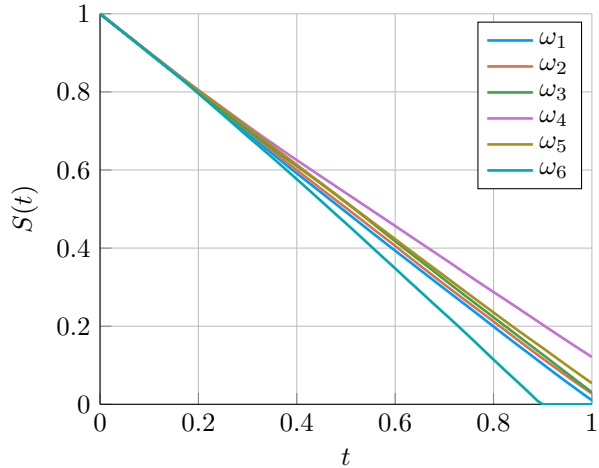}
  \caption{Sample paths of $S^\alpha(t)$ following the GBM
    model~\eqref{eq:gbm_sde}, with $q(\alpha(t))=1$, and GBM volatility
    $\sigma = 0.1$. The paths are  more regular than that of the
    Brownian noise model shown in \Cref{fig:brownian_paths}.
  }\label{fig:sde_gbm_const_x}
\end{figure}

We therefore argue that it is generally inappropriate to model the uncertainty
in future demand with increments of a Brownian motion. Instead, the
uncertainty in demand could be modelled more realistically by
describing the evolution of parameters in the demand function.
This article will focus on multiplicative demand
uncertainty, which may arise from estimates of seasonality or
long-term trends in demand. As such, we introduce a multiplicative
parameter $g\geq 0$, and consider the demand function given by
$(a,g)\mapsto q(a)g$. Assume that the parameter is not directly
affected by price, but  only represents exogenous information outside the
retailer's control.
For our purposes, we will model the multiplicative parameter as
a geometric Brownian motion (GBM), $G(t) = \exp\left(
  -\frac{\sigma^2}{2}t +\sigma W(t)\right)$, with volatility coefficient
$\sigma\geq 0$. The GBM approach is also considered
by~\cite{xu2006monopolistic} for a different dynamic pricing problem from
that presented in this article.
We restrict ourselves to a GBM with no
drift to focus on the impact of the uncertainty, rather than modelling
time-dependent behaviour such as seasonality.
For any $\Delta t\geq 0$, some relevant properties of $G(t)$ are as follows.
\begin{align}
  G(0)&=1,\\
  \mathbb E[G(t+\Delta t)\mid G(t)]&=G(t),
                                   &\text{martingale property,}\\
  \mbox{Var}[G(t+\Delta t)\mid G(t)]&={G(t)}^2(e^{\sigma^2 \Delta
                                      t}-1), & \text{increasing variance.}
\end{align}
With this model, we expect future demand $q(a)G(t)$ at a given price $a\geq
0$ to be the current experienced demand, but with decreasing certainty
the further ahead we forecast.
The SDE governing the system, started at
$S(0)>0$, $G(0)=1$, is
\begin{align}
  \begin{split}\label{eq:gbm_sde}
    \mathrm{d}G(t)&=\sigma G(t)\,\mathrm{d}W(t),\\
    \mathrm{d}S^\alpha(t)&=-q(\alpha(t))G(t)\,\mathrm{d}t,\qquad\text{stopped at zero}.
  \end{split}
\end{align}
The sample paths of $S^\alpha(t)$ following the GBM model, as shown in
\Cref{fig:sde_gbm_const_x}, are more regular
than of the Brownian noise model.

\subsection{Non-dimensionalized system}\label{subsec:nondimensionalization}
In order to capture similarities between different pricing decisions, irrespective of units
such as a particular currency, it is helpful to work in a dimensionless system.
We thus non-dimensionalize the model, which also helps to reduce the
number of parameters in the decision problem.
The units at play in our system are the time and the product price,
for example measured in weeks and \pounds.
If $s,a$, and $t$ denotes unscaled quantities of stock, price, and
time respectively, we rescale them with dimensionless hat-quantities
\begin{align}
  \hat s&=\frac{s}{S(0)},&\hat t &=\frac{t}{T},
  &\hat a &= \frac{a}{\bar a}.
\end{align}
Here, $S(0)\gg 1$ is the initial quantity of stock, $T$ is the time-horizon
for the pricing problem, and $\bar a$ is some reference price chosen
to make typical prices $\hat a$ continuous and of order one.
In order to write down a dimensionless formulation of the system's
SDE~\eqref{eq:gbm_sde}, we need to work with the expected demand
function and volatility parameter defined as
\begin{align}
  \hat q(\hat a)&=\frac{T}{S(0)}q(\hat a \bar a),
  &\hat \sigma &= \sqrt{T}\sigma.
\end{align}
Thus, for a given pricing policy $\hat \alpha(\hat t)$, the
system starts at $\hat S(0),\hat G(0)=1$, and evolves according to
\begin{align}
  \begin{split}\label{eq:gbm_sde_nondim}
    d\hat G(\hat t)&=\hat \sigma \hat G(\hat t)\,\mathrm{d}W(\hat t),\\
    d\hat S^{\hat \alpha}(t)&=-\hat q(\hat \alpha(\hat t))\hat G(\hat
    t)\,\mathrm{d}\hat t,\qquad\text{stopped at zero}.
  \end{split}
\end{align}
This scaling means that we work solely with stock and time on the unit interval, $\hat S^{\hat \alpha}(\hat t),\hat t\in[0,1]$.
In the remainder of the article, we will omit the hats, and assume all
quantities are dimensionless.

\section{Formulation of the decision
  problem}\label{sec:decision_formulation}
We restate the retailer's objective which
guides the choice of pricing policy: Over some decision horizon, continuously
adjust the price of a given product in order to maximize
the total profit, generated from sales revenue minus the cost of handling
unsold items at the terminal time. We formulate this mathematically as
a stochastic optimal control problem, for which the optimal pricing
policy can be found by solving the associated HJB equation.
In \Cref{subsec:param_estim} we address the real-world restrictions of
the continuous-time assumption.

Let $C>0$ denote the handling cost per unit of stock at the terminal
time. Define $T_h$ to be the hitting time $T_h=\min \{1, T_0\}$, where
$T_0=\inf\{t\geq 0\mid S^\alpha(t)=0\}$.
The total profit accrued from a pricing policy $\alpha(t)$ is then
\begin{equation}\label{eq:profit_expression}
  P(\alpha) = \int_0^{T_h}\alpha(u)q(\alpha(u))G(u)\,\mathrm{d}u - CS^\alpha(1).
\end{equation}
Note that we focus on time horizons where we believe discounting
future cash is negligible.
This profit is a random variable that depends on the event $\omega$
or, equivalently, the path of the
Brownian motion $W(t)$. 
In the context of this article, we assume a retailer that
wants to maximize the expected value of profit, $\mathbb
E[P(\alpha)]$.

Further, we restrict the product price to be in some closed  interval
$A\subset\mathbb R_{\geq 0}$.
In addition, we will only look for Markovian pricing policies in an
admissible set $\mathcal A$. We say that $\alpha(t)$ is Markovian if
there is a function of the form $a(t,s,g)$ such that for each event $\omega$,
\begin{equation}
  \alpha(t)(\omega) = a(t,S^\alpha(t)(\omega), G(t)(\omega)) \in A.
\end{equation}
Thus, we seek pricing policies that set the price at time $t$, based
on the knowledge of the state at that time. We also assume that
$\mathcal A$ only contains pricing policies such that the integral
in~\eqref{eq:profit_expression} exists and
$P(\alpha)$ is integrable. We can now state the mathematical problem we
seek to solve in the remainder of the article.
\begin{mydef}[Pricing problem]
  In order to maximize the retailer's expected profit, find a solution to
  the stochastic optimal control problem
  \begin{equation}\label{eq:pricing_problem}
    \max_{\alpha\in \mathcal A}\,\mathbb E[P(\alpha)].
  \end{equation}
\end{mydef}
Our chosen strategy for finding the corresponding pricing function $a^*(t,s,g)$
of a maximizer $\alpha^*\in\mathcal A$ is to solve the associated HJB
equation for the pricing problem.
For a detailed explanation of the theory behind stochastic optimal control and
HJB equations, see, for example,~\cite{pham2009continuous}.
The HJB equation is a nonlinear PDE, where the solution describes the
value of being in a particular state. We will not worry about
uniqueness of solutions to the HJB equation in this article. From the
value function defined below, one can
calculate the optimal pricing function $a^*(t,s,g)$.

If at time $t$, we know the value of $S(t)$ and $G(t)$, then the value
function represents the expected value of applying the optimal pricing
policy for the remainder of the pricing period.
We define the value function $v(t,s,g)$ for $t,s\in[0,1]$ and $g\geq 0$ by
\begin{align}
  v(t,s,g) &=
             \max_{\alpha \in \mathcal A}\,
             \mathbb E_{t}\left[
             \int_t^{T_h}\alpha(u)q(\alpha(u))G(u)\,\mathrm{d}u-CS^\alpha(1)
             \right],\label{eq:def_value_interior}\\
  v(t,0,g)&= 0,\\
  v(t,s,0)&=-Cs.
\end{align}
The subscript on the expectation denotes that we condition on
$S^\alpha(t)=s$ and $G(t)=g$.
In the limit $t\to 1$, we see from~\eqref{eq:def_value_interior} that $v(1,s,g)=-Cs$.
In the usual way, see for example~\cite{pham2009continuous}, the HJB equation for $v$ is
\begin{equation}\label{eq:hjb_interior}
  v_t(t,s,g)+\frac{\sigma^2}{2} g^2v_{gg}(t,s,g)
  +g\max_{a\in A}\{(a-v_s(t,s,g))q(a)\} = 0, \quad
  (t,s,g)\in{(0,1)}^2\times\mathbb R_{>0}.
\end{equation}
The subscripts on $v$ denote partial derivatives with respect to the
given argument. Together with the boundary and terminal conditions on
$v$, this constitutes the HJB equation for the pricing problem.
If we know $v$, the optimal pricing function can be calculated from
the univariate optimization problem
\begin{equation}\label{eq:pricing_function_general}
  a^*(t,s,g) = \argmax_{a\in A}\{(a-v_s(t,s,g))q(a)\}.
\end{equation}
In \Cref{sec:deterministic_hjb,sec:stochastic_hjb}, solutions for the
pricing problem are found via the HJB equations with linear and
exponential demand functions.
\begin{myremark}\label{rem:viscosity}
  The value function is not necessarily sufficiently smooth to satisfy
  the HJB equation~\eqref{eq:hjb_interior} in the classical sense.
  This is indeed the case for some examples in this article when $\sigma=0$,
  and we must therefore consider
  the solutions in the viscosity sense. See~\cite{pham2009continuous}
  for a description of viscosity solutions to HJB equations.
\end{myremark}

\subsection{Parameter estimation}\label{subsec:param_estim}
In a real retail application we cannot update the price in continuous
time, and it is not possible to infer $G(t)$ exactly.
To represent real world conditions, we assume that
the price is piecewise constant and updated frequently at fixed time points
$t_0<t_1<\cdots<T_h$.
Further, we assume that stock levels are only observed at these time
points. When computing an optimal pricing strategy
we will still consider the continuous-time
function~\eqref{eq:pricing_function_general} that can be computed from
the HJB equation as
the optimal pricing function for the problem. Such a continuous-time
approximation to inherently discrete systems is common, for example
in pricing of options using the Black-Scholes equations
\citep{black1973pricing}.
In order to use a pricing function $a(t_k,s,g)$ we must estimate
$G(t_k)$.
Due to the Markov properties of $G(t)$ and $S^{\alpha}(t)$, information
about the process for $t<t_{k-1}$ is not needed, and we can estimate
$G(t_k)$ based on $\alpha(t_{k-1})$, $S^\alpha(t_{k-1})$, and $S^\alpha(t_k)$.
Let $k = 1$ and leave out the superscript $\alpha$ of $S^\alpha(t)$ for the
remainder of this section.
By assumption, the price has been constant, $\alpha(u)=a_0$, for the time
period $u\in[t_0,t_1)$. Say $S^\alpha(t_1)>0$, and that we wish to update the price at time
$t_1$.
From~\eqref{eq:gbm_sde_nondim}, the SDE describing our system, we have
\begin{equation}\label{eq:s_integrated}
  S(t_1)=S(t_0)-q(a_0)\int_{t_0}^{t_1}G(u)\,\mathrm{d}u.
\end{equation}
In the numerical examples in this article, we estimate $G(t_1)$ with
$\hat G(t_1)=\frac{S(t_0)-S(t_1)}{q(a_0)(t_1-t_0)}$.
We now discuss the derivation and properties of this estimator.

Define $B(t) =  \exp(-\sigma^2 t/2+\sigma\sqrt{t}Z_0)$, where $Z_0\sim
\mathcal{N}(0,1)$. The evolution of $G(t)$ is known in closed form,
which gives
$G(t_1)=G(t_0)B(\Delta t)$, with $\Delta t = t_1-t_0$.
From~\eqref{eq:s_integrated} it follows that
\begin{align}
  G(t_1) = \frac{S(t_0)-S(t_1)}{q(a_0)}\frac{B(\Delta
  t)}{\int_0^{\Delta t}B(u)\,\mathrm{d}u}.
\end{align}

\begin{figure}[htbp]
  \centering
  \includegraphics[width=0.5\textwidth,axisratio=1.55]{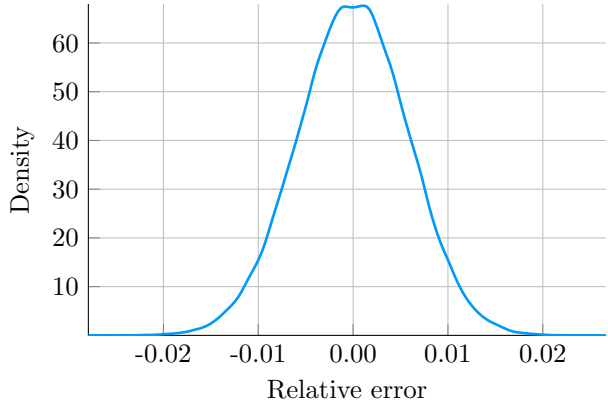}
  \caption{Distribution of the relative error in estimating $G(t)$ with
    $\sigma=0.1$ and $\Delta t= 0.01$, based on \num{100000} samples.}\label{fig:gbm_estimate}
\end{figure}

So long as $\Delta t=t_1-t_0$, and the variance of $B(u)$, are
sufficiently small, we may use the approximation
$\int_{0}^{\Delta{t}}B(u)\,\mathrm{d}u\approx  \frac{\Delta
  t}{2}(B(0)+B(\Delta t))=\frac{\Delta
  t}{2}(1+B(\Delta t))$.
Hence, the conditional distribution of $G(t_1)$ is approximated as
\begin{align}
  G(t_1)\mid (a_0,S(t_0),S(t_1))
  &\approx
    \frac{S(t_0)-S(t_1)}{q(a_0)}\frac{2}{\Delta t}\frac{B(\Delta t)}{1+B(\Delta t)}\\
  &\approx \frac{S(t_0)-S(t_1)}{q(a_0)}\frac{1+B(\Delta t)}{2\Delta t}.
\end{align}
The second approximate equality comes from the Taylor
expansion $\frac{x}{1+x}\approx \frac{1}{4}(1+x)$ about $x=1$.
This gives us the following expressions for the first two moments
of $G(t_1)$:
\begin{align}\label{eq:conditional_mean_g}
  \mathbb E[G(t_1)\mid a_0,S(t_0),S(t_1)]
  &\approx \frac{S(t_0)-S(t_1)}{ q(a_0)\Delta t},\\
  \mbox{Var}[G(t_1)\mid a_0,S(t_0),S(t_1)]
  &\approx \frac{1}{4}{\left(\frac{S(t_0)-S(t_1)}{
    q(a_0)\Delta t}\right)}^2\left( e^{\sigma^2\Delta t} -1 \right).
\end{align}
The conditional expectation in~\eqref{eq:conditional_mean_g} is equal
to our estimator $\hat G(t)$.
\Cref{fig:gbm_estimate} shows the distribution of the relative
difference $1-\hat G(t)/G(t)$ between
the estimate and the true $G(t)$, for
the parameters used in \Cref{sec:stochastic_hjb}. The relative
estimator error when $\sigma=0.1$ and $\Delta t=0.01$
is typically within \SI{1}{\percent}.

\section{Solution to the deterministic
  system}\label{sec:deterministic_hjb}
In this section, we provide solutions to the pricing problem in the
deterministic case, $\sigma = 0$, for families of linear and exponential
demand functions $q_l(a)$ and $q_e(a)$ respectively.
\begin{align}
  q_l(a)&=q_1-q_2a,&\text{for } q_1,q_2>0,\\
  q_e(a)&=q_1e^{-q_2a},&\text{for } q_1,q_2>0.
\end{align}
These demand functions are often used in the literature. For a discussion
about their properties and usage in modelling demand, see~\cite[Ch.~7]{talluri2006theory}.
With both demand functions, the optimal pricing function policy is to hold the
price constant over the pricing horizon. They translate to the
following intuition: For small amounts of stock, relative to demand,
sell at the highest price so that all stock is depleted by the
terminal time. As the amount of stock increases, the retailer should
decrease the price, until it arrives at the price that maximizes the
profits balancing revenue and cost of unsold stock at the terminal
time. The ansatz pricing policies referred to below come from the
discrete-time pricing problem considered by~\cite{riseth2017comparison}.

\subsection{Linear demand function}
When $q(a)=q_1-q_2a$ for $q_1,q_2>0$, the maximization in the HJB
equation can be solved in closed form.
To reduce the number of parameters, we can rescale the price per
product with $q_2$ by setting $q_2a$ to $a$.
Then $q(a)=q_1-a$, and we set the pricing interval to $A=[0,q_1]$ so that the demand
function is non-negative.
Now, use the ansatz that $a^D(t,s,g)$ sets the price so that
$\alpha^D(t)$ is
constant for the remainder of the pricing period.
From the expression of the value function
in~\eqref{eq:def_value_interior}, the ansatz pricing policy must
also be given by
\begin{align}
  a^D(t,s,g)
  &=\argmax_{a\in A}\left\{\int_t^1aq(a)g\,\mathrm{d}u -
    C\left(s-\int_t^1q(a)g\,\mathrm{d}u\right) \,\big|\,  s \geq \int_t^1q(a)g\,\mathrm{d}u\right\}\\
  &=\argmax_{a\in A}\left\{(1-t)(a+C)q(a)g-Cs \mid  s \geq
    (1-t)q(a)g\right\}.
\end{align}
The maximizer therefore satisfies the equality constraints that
$s$ equals $(1-t)gq(a)$, is zero, or is in the interior of the feasible set,
given by $A$.
For $q(a)=q_1-a$, we can verify that this implies
\begin{equation}\label{eq:astar_linear}
  a^D(t,s,g)=\begin{cases}
    q_1-\frac{s}{(1-t)g}, &\text{if } 0\leq s\leq
    (1-t)g\min\{q_1,\frac{1}{2}(q_1 +C)\},\\
    \max(0,q_1-C)/2,&\text{otherwise.}
  \end{cases}
\end{equation}
It follows that $T_h=1$ when following the optimal pricing strategy.
This is an intuitive result, because if $T_h<1$, one can increase the
price until $T_h=1$ and $S^\alpha(1)=0$, which earns extra revenue at no extra cost.
The pricing policy suggested by the
deterministic assumption provides the following, obvious, heuristic:
First, find the price that maximizes profits, ignoring inventory
constraints. Second, if sales forecasts suggest you will deplete stock
before the end of the time horizon at this price, increase it
accordingly. This is consistent with the heuristic proposed
by~\cite{schlosser2015dynamic1},~\cite{schlosser2015dynamic2}, which
he finds by considering a deterministic continuum
approximation to a Poisson process demand model.

\begin{myex}\label{ex:acecplot}
  To demonstrate what form the pricing function may take,
  \Cref{fig:hjbsol_a_cec} shows a plot of $a^D(t,s,1.0)$ for a given set of
  parameters. We have chosen the parameters so that the pricing
  problem starts at the most interesting point: at the kink separating
  the regimes where all the stock is sold out and where it is not.
  This corresponds to a combination of parameters such that
  $\min(q_1,\frac{1}{2}(q_1+C))=1$.

  \begin{figure}[htbp]
    \centering
    \includegraphics[width=0.6\textwidth,axisratio=1.15]{hjbsol_a_cec}
    \caption*{$a^D(t,s,1.0)$}
    \caption{Example optimal, deterministic, pricing function $a^D(t,s,1.0)$
      from~\eqref{eq:astar_linear}, as a function of time and stock.
      Notice how sensitive $a^D(t,s,1.0)$ is to changes in $s$, for large
      $t$.
      The parameters used are $q_1=3/2$, and $C=1/2$.
      In \Cref{fig:hjbsol_a} we compare this function to the optimal price
      when $\sigma=0.1$.
    }\label{fig:hjbsol_a_cec}
  \end{figure}
\end{myex}

To verify that $a^D(t,s,g)$ given by~\eqref{eq:astar_linear} indeed is
the solution to the deterministic
pricing problem, we show that it satisfies the HJB equation.
On the interior of the $(t,s,g)$-domain, $a^D(t,s,g)$ must solve
\begin{equation}
  \max_{a\in A} \{(a-v_s(t,s,g))(q_1-a)\}.
\end{equation}
Let $\mathcal P_A$ denote the projection of the real line to $A$. The
objective is concave, and hence the maximizer is
\begin{equation}\label{eq:astar_hjb_linear}
  a^D(t,s,g) = \mathcal P_A\left[\frac{q_1+v_s}{2}\right].
\end{equation}
Define $\Gamma$ to be the boundary of the terminal time problem,
that is when $t=1$, $s=0$, or $g=0$.
Let us assume that $C< q_1$, then the
deterministic-case HJB equation for $v$ is
\begin{align}\label{eq:hjb_linear}
  v_t(t,s,g)+
  \frac{g}{4}{( q_1-v_s(t,s,g))}^2
  &=0,&(t,s,g)&\in{(0,1)}^2\times \mathbb R_{>0},\\
  v(t,s,g) &= -Cs,&(t,s,g)&\in \Gamma.
\end{align}
If $C\geq q_1$, then there are regions where $a^D(t,s,g)=0$. In
particular, this means that for $t,s,g$ such that
$s>(1-t)gq_1$, $v$ must satisfy
\begin{equation}\label{eq:hjb_linear_Clarge}
  v_t(t,s,g) - gq_1v_s(t,s,g)=0.
\end{equation}

From the two expressions for $a^D(t,s,g)$ in~\eqref{eq:astar_linear}
and~\eqref{eq:astar_hjb_linear}, the ansatz implies that the value
function must satisfy
\begin{align}\label{eq:v_func_linear_det}
  v^D(t,s,g)
  &=\begin{cases}
    q_1s-\frac{s^2}{(1-t)g},
    &\text{if } 0\leq s\leq
    (1-t)g\min\{q_1,\frac{1}{2}(q_1 +C)\},\\
    -Cs + V(C)(1-t)g,
    &\text{otherwise.}
  \end{cases}\\
  V(C)
  &= \begin{cases}
    {\left[\frac{1}{2}(q_1+C)  \right]}^2,&C<
    q_1,\\
    q_1C,&C\geq q_1.
  \end{cases}
\end{align}
This $v$ does indeed satisfy the deterministic HJB equation given
by~\eqref{eq:hjb_linear}-\eqref{eq:hjb_linear_Clarge}, and hence we
can conclude that $a^D(t,s,g)$ is an optimal pricing function. Note that $v$ is
not smooth for all parameter combinations, and is therefore considered
a solution in the viscosity sense as noted in Remark~\ref{rem:viscosity}.

\subsection{Exponential demand function}
With the same ansatz that was used for the linear demand function, we
can also find the optimal, deterministic, pricing function for
exponential demand $q(a)=q_1e^{-q_2a}$.
Indeed, this is true for any demand function for which a closed form
solution exists for $\max_{a\in A}\{(a+C)q(a))\}$ and $(1-t)gq(a)=s$.
As with the linear demand, we can eliminate the parameter $q_2$ so that
$q(a)=q_1e^{-a}$, given by replacing $q_2a$ by $a$.
In the exponential demand case, with $A=[0,\infty)$, the ansatz gives
us the optimal pricing function
\begin{equation}\label{eq:afun_exp_cec}
  a^D(t,s,g)=\begin{cases}
    \log{\frac{q_1g(1-t)}{s}},
    &\text{if } 0\leq\frac{s}{q_1g(1-t)}\leq e^{C-1},\\
    \max\left(0,1-C\right),&\text{otherwise}.
  \end{cases}
\end{equation}
For completeness, we state the HJB equation for the exponential demand
case when $C<1$, and provide the solution so that the optimality of~\eqref{eq:afun_exp_cec} can
be verified. The maximizer of $\max_{a\in A}\{(a-v_s)q(a)\}$ in the
HJB equation is $a^D =  1+v_s$. Thus, the value function
must satisfy the HJB equation
\begin{align}
  v_t(t,s,g)
  +gq_1e^{-1-v_s(t,s,g)}
  &=0,&(t,s,g)&\in{(0,1)}^2\times \mathbb R_{>0},\\
  v(t,s,g) &= -Cs,&(t,s,g)&\in \Gamma.
\end{align}
The viscosity solution of the HJB equation, acquired from the ansatz $a^D(t,s,g)$
in~\eqref{eq:afun_exp_cec}, is
\begin{align}
  v^D(t,s,g)
  &=\begin{cases}
    s \log \frac{q_1g(1-t)}{s},
    &\text{if } 0\leq \frac{s}{q_1g(1-t)}\leq e^{C-1},\\
    -Cs+gV(C)(1-t),&\text{otherwise},
  \end{cases}\\
  V(C)&=q_1e^{C-1}.
\end{align}

\section{Impact of uncertainty}\label{sec:stochastic_hjb}
We now discuss to what degree multiplicative uncertainty changes
our policy. In this section, we solve the HJB equation numerically with the linear
demand function $q(a)=q_1-a$ defined on $A=[0,q_1]$, and compare the
resulting pricing policy to the
deterministic-system policy from the previous section.
With the diffusion term in the HJB equation, one can
expect the kink in the deterministic pricing function to smooth
out. It turns out that the difference between an optimal
policy and a heuristic policy based on the solution
to the deterministic system is at most
$\mathcal{O}(\sigma \sqrt{1-t})$. Further, numerical tests indicate that the
closed-form pricing functions found in the previous section
perform sufficiently well in most situations.

We assume in the following that $(q_1+v_s)\in[0,2]$ for
$t\in[0,1)$, $s,g>0$, so
that the pricing function satisfies $a(t,s,g) = (q_1+v_s(t,s,g))/2$ as
given by~\eqref{eq:astar_hjb_linear}.
Using~\eqref{eq:hjb_interior} and~\eqref{eq:hjb_linear} it then
follows that $v$ should satisfy
\begin{align}\label{eq:hjb_impact_uncertainty}
  v_t(t,s,g)+\frac{\sigma^2}{2} g^2v_{gg}(t,s,g)
  +\frac{g}{4}{(q_1-v_s)}^2
  &= 0,&(t,s,g)&\in{(0,1)}^2\times \mathbb R_{>0},\\
  v(t,s,g) &= -Cs,&(t,s,g)&\in \Gamma.
\end{align}

The numerical solution to the HJB equation is solved with the
following procedure: (i)~Reformulate the PDE with the similarity
transformation $\xi = s/g$ and $v = g\phi$. (ii)~Truncate the boundary
for $\xi\to\infty$ and set an asymptotic Dirichlet boundary condition
based on the deterministic-system solution. (iii)~Approximate
the PDE for $\phi(t,\xi)$ with central finite differences and the
\texttt{Tsit5} time stepping procedure in
\texttt{DifferentialEquations.jl}
\citep{rackauckas2017differentialequations}, implemented in the Julia
programming language \citep{bezanson2017julia}.
We denote the computed pricing function and pricing
policy by $a^B(t,s,g)$ and $\alpha^B(t)$ respectively.

\begin{myex}\label{ex:function_impact_uncertainty}
  Let us consider the particular example system used in
  Example~\ref{ex:acecplot}. 
  That is, a linear demand function $q(a)=q_1-a$, with $q_1=3/2$
  and $C=1/2$.
  We set the volatility level of $G(t)$ to $\sigma=0.1$, which
  corresponds to a true demand near the terminal time within \SI{20}{\percent} of the
  expected demand $q(a)$, with probability $0.95$.
  \Cref{fig:hjbsol_a} shows the
  optimal pricing function for $g=1$, and a plot of the difference
  $a^B(t,s,g)-a^D(t,s,g)$. The only visible difference is along the kink
  line, $s=g(1-t)$, where $a^B$ smooths out the transition between
  the two regions, and hence sells the product at a slightly higher
  price.
  \begin{figure}[htb]
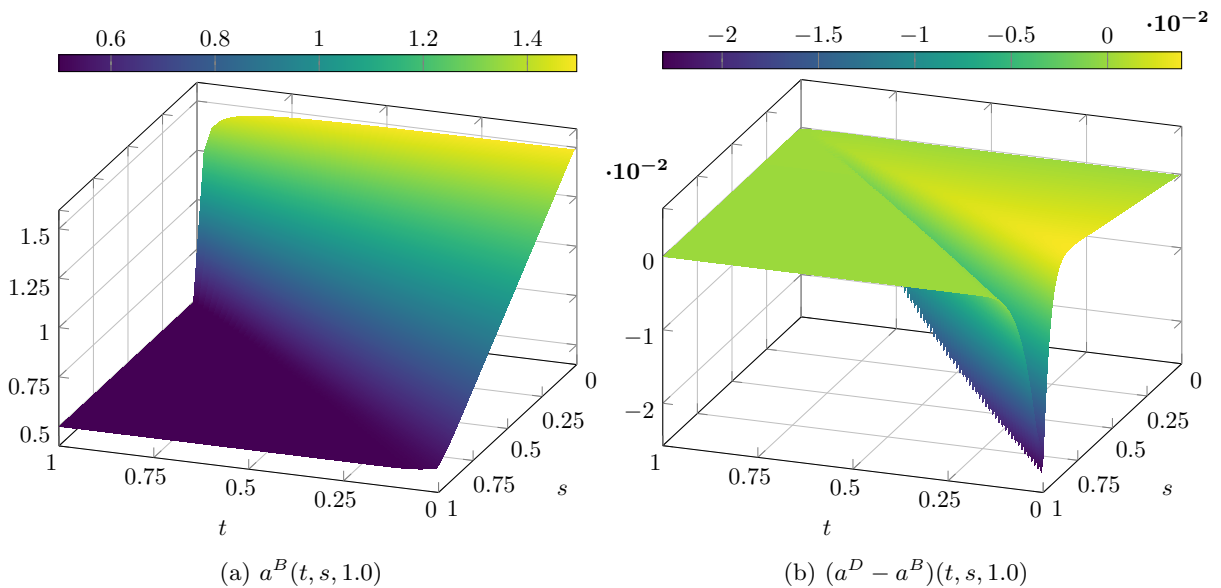

    \centering
    \begin{subfigure}[b]{0.5\textwidth}
      \includegraphics[width=\textwidth]{hjbsol_a_hjb}
      \caption{$a^B(t,s,1.0)$}
    \end{subfigure}%
    \begin{subfigure}[b]{0.5\textwidth}
      \includegraphics[width=\textwidth]{hjbsol_a_cec_hjb}
      \caption{$(a^D-a^B)(t,s,1.0)$}
    \end{subfigure}%
    \caption{The optimal pricing function ($a^B(t,s,g)$, left) smooths out the
      kink as compared to the deterministic heuristic ($a^D(t,s,g)$,
      \Cref{fig:hjbsol_a_cec}).
      The right plot shows the impact of uncertainty on the optimal
      pricing function: (i)~When we do not expect to sell out of the
      product, the prices are the same. (ii)~When we expect to sell out
      of the product, the deterministic heuristic takes a slightly larger
      price.
      (iii)~In the transition between the two regions, uncertainty
      increases the optimal price.
    }\label{fig:hjbsol_a}
  \end{figure}
\end{myex}

\begin{figure}[hbt]
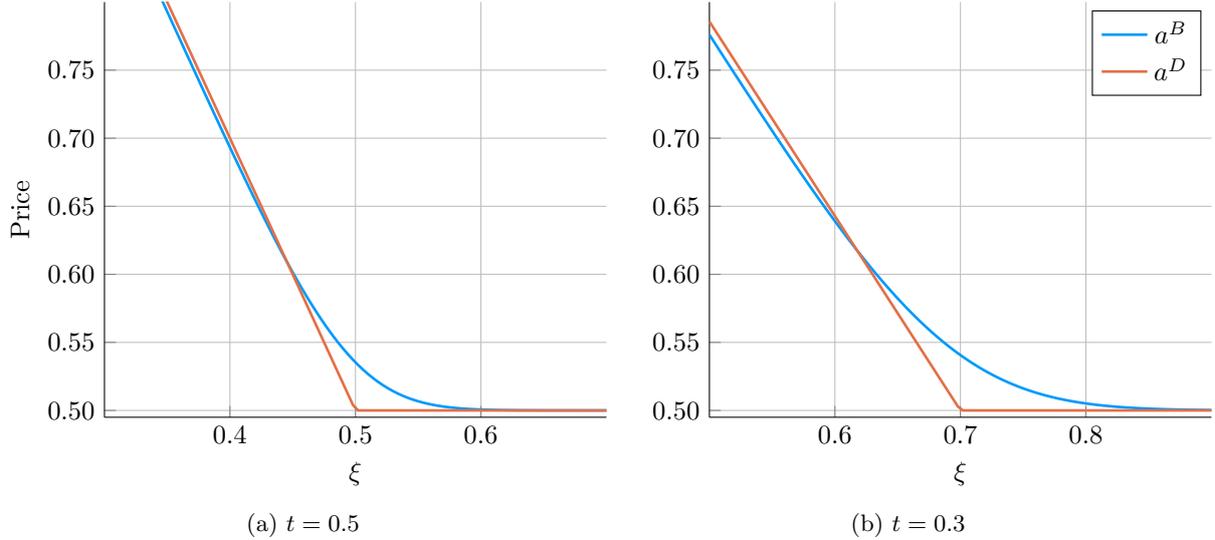

  \centering
  \begin{subfigure}[b]{0.5\textwidth}
    \includegraphics[width=\textwidth,axisratio=1.2]{hjb_cec_zoom_02_05}
    \caption{$t = 0.5$}
  \end{subfigure}%
  \begin{subfigure}[b]{0.5\textwidth}
    \includegraphics[width=\textwidth,axisratio=1.2]{hjb_cec_zoom_02_07}
    \caption{$t=0.3$}
  \end{subfigure}
  \caption{A zoomed in comparison of $a^B(t,s,g)$ and $a^D(t,s,g)$
    in the transformed coordinate $\xi=s/g$.
    The difference between $a^B$ and $a^D$ is of order
    $\sigma\sqrt{1-t}$ at the kink, and
    $-\sigma^2\xi/2$ as we move to the left in the
    plots.
    This example uses $q_1=3/2$, $C=1/2$, and $\sigma=0.2$.
  }\label{fig:hjb_cec_zoom}
\end{figure}
\Cref{fig:hjbsol_a} indicates that one can do an asymptotic analysis
of the impact of $0<\sigma \ll 1$ on the pricing function. The
appendix provides details of the analysis. We summarise the
results here, and refer to \Cref{fig:hjb_cec_zoom} for a visualisation
of the differences between $a^B$ and $a^D$.
\begin{myresult}\label{res:linear_asymptotics}
  Define $\beta=\min(q_1,(q_1+C)/2)$.
  There is an inner layer around the surface $\frac{s}{g}=(1-t)\beta$ which
  smooths out the kink of the solution $a^D(t,s,g)$ that arises when
  $\sigma=0$. At the kink, the first-order correction in the pricing
  function is of order $\mathcal{O}(\sigma\sqrt{1-t})$.
  The width of the layer is of the order
  $\mathcal{O}(\sigma{(1-t)}^{3/2})$, which connects the two pricing
  regimes identified in the deterministic case.
  As we move from the inner layer to larger values of
  $\xi:=s/g>(1-t)\beta$, the solution
  tends to the value $\delta=\max(0,q_1-C)/2$, which coincides with
  $a^D$. As we move from the inner layer to smaller values of $\xi$,
  the leading order solution of the inner layer coincides with $a^D$.
  Further, there is a second order correction in the region
  $\xi<(1-t)\beta$ equal to $-\sigma^2 \xi/2$.
\end{myresult}

The first and second order corrections reflect the insights one can
arrive at when taking into account the deviations of actual future
demand from expected demand.
At the interface where we expect to sell all the inventory at the
optimal lower-bound, ``infinite-inventory'' price $\delta$, taking
into account the possibility of higher future demand means that we can
increase the price.
When the inventory is so low that we expect to sell it all at some
price higher than $\delta$, taking into account the possibility
of lower future demand means that we should price the product lower than in the
deterministic-case in order to
reduce the probability of having excess inventory at the terminal time.

The asymptotic results together with \Cref{fig:hjbsol_a} indicate that one can expect a price
decrease over time when following the optimal pricing policy
$\alpha^B(t)$.
The numerical investigation in the next example verifies this, but highlights
that there are negligible gains in total profit from pricing according to
$\alpha^B(t)$ rather than the deterministic-case
pricing policy $\alpha^D(t)$.

\begin{figure}[htb]
    \centering
    \includegraphics{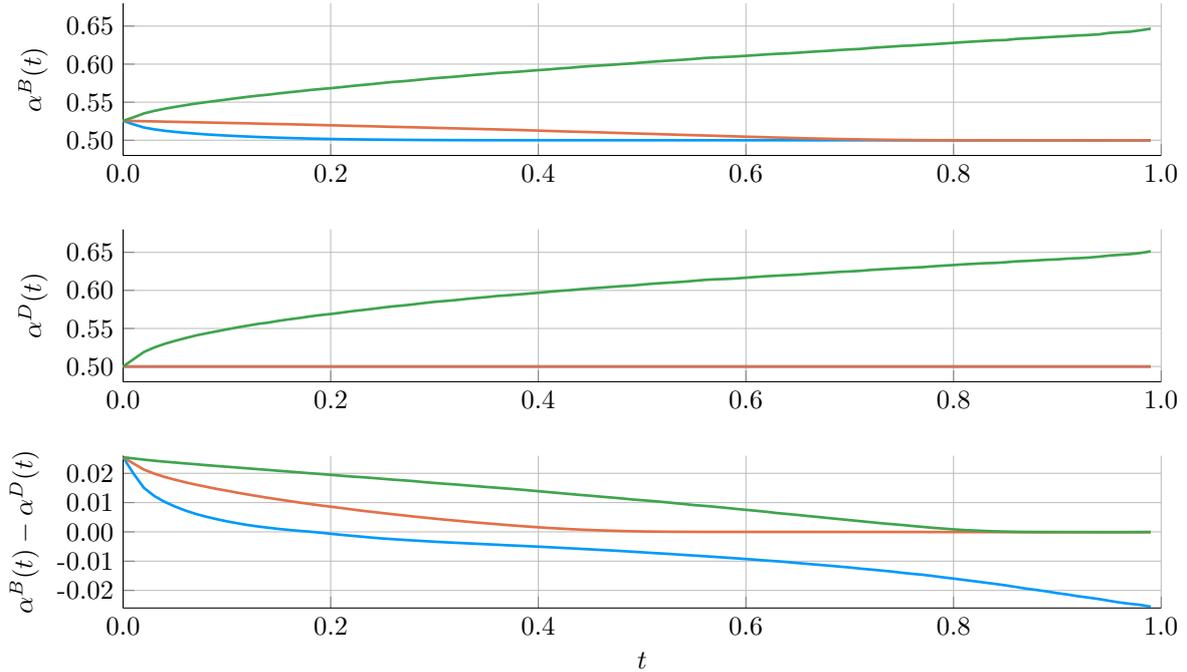}
    \caption{Statistics of the price paths from
      Example~\ref{ex:impact_uncertainty}, 
      showing the $0.05,0.5$, and
      $0.95$ quantiles of $\alpha(t)$. From top to bottom, the optimal
      policy, the deterministic heuristic, and their difference.
      The optimal policy starts slightly higher, and decreases over
      time.
      Note that the $0.05$ and $0.5$ quantiles lie on top of each
      other for $\alpha^D(t)$, so prices are likely to stay constant
      at $1/2$.
    }\label{fig:gbm_stats_a}
  \end{figure}
\begin{myex}\label{ex:impact_uncertainty}
  Let us simulate the system from
  Example~\ref{ex:function_impact_uncertainty}
  with the numerical approximation to $\alpha^B(t)$ and compare it to $\alpha^D(t)$.
  We draw \num{100000} sample paths from the underlying Brownian
  motion $W(t)$, and set the price to be constant on intervals of size $\Delta t = 0.01$
  using a policy $\alpha(t)$. We use the estimator described in
  \Cref{subsec:param_estim} as an approximation to $G(t)$.
  The simulations are run with both $\alpha^B(t)$ and $\alpha^D(t)$, and
  statistics of their paths are shown in \Cref{fig:gbm_stats_a}.
  The prices $\alpha^B(t)$ start slightly higher than $\alpha^D(t)$, but
  will then over time decrease, on average, towards the lower bound
  $\delta=1/2$. We also see that the deterministic heuristic is less
  anticipative, and will begin increasing the prices compared to the
  optimal policy after $t=0.2$.

  The measure of interest, however, is how much profit the different
  policies make. Recall that the profit of following a policy
  $\alpha(t)$ is the random variable
  \begin{equation}
    P(\alpha)=\int_0^{T_h}q(\alpha(u))\alpha(u)G(u)\,\mathrm{d}u-CS^\alpha(1).
  \end{equation}
  From simulations we estimate the distribution of $P(\alpha)$ for the
  optimal and deterministic policies, and compare their performance.
  The improvement is negligible, as we see from the relative statistics
  \begin{align}
    \mathbb E\left[1-P(\alpha^D)/P(\alpha^B)\right]
    &\approx 2\times 10^{-4},
    &\mbox{std}\left[1-P(\alpha^D)/P(\alpha^B)\right]
    &\approx 6\times 10^{-4},\\
    \mbox{Median}\left[1-P(\alpha^D)/P(\alpha^B)\right]
    &\approx -1\times 10^{-4}.
  \end{align}
  The calculated optimal pricing policy results in \SI{0.02}{\percent} higher profits
  than the heuristic, on average. It even results in lower profits than the
  heuristic more than \SI{50}{\percent} of the realizations.
  \Cref{fig:gbm_histogram_a_cec_hjb} shows a histogram that approximates
  the distribution of the relative loss from using $\alpha^D(t)$
  over the optimal pricing policy. The differences between the two are
  small, but the distribution is non-symmetric:
  The heuristic $\alpha^D(t)$ results in slightly larger profit for more
  than half of
  the realizations, at the expense of performing  worse for the remaining realizations.
  \begin{figure}[htb]
    \centering
    \includegraphics[width=0.6\textwidth,axisratio=1.25]{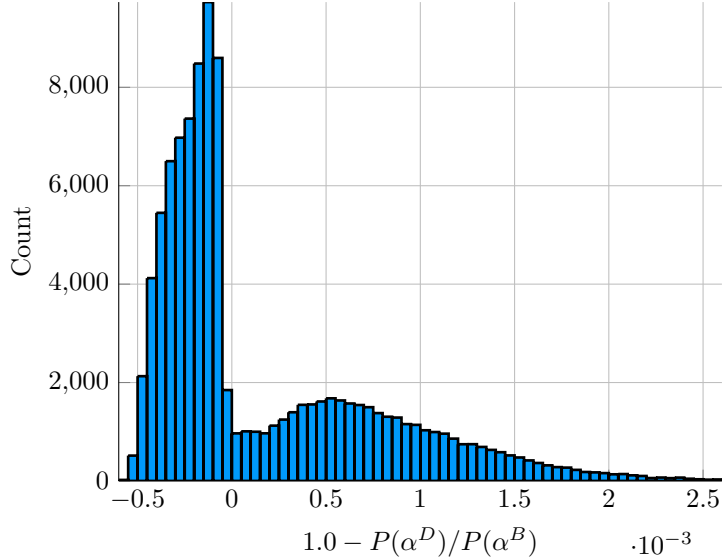}
    \caption{The histogram shows the distribution of the relative
      profit by using the deterministic pricing heuristic $\alpha^D(t)$
      to using the optimal policy $\alpha^B(t)$, as described in
      Example~\ref{ex:impact_uncertainty}. 
      Positive values
      correspond to realizations of $W(t)$ where
      $\alpha^B(t)$ is better.
      The shape of the distribution is similar to the discrete-time pricing problem
      studied by~\cite{riseth2017comparison}.
    }\label{fig:gbm_histogram_a_cec_hjb}
  \end{figure}
\end{myex}

In Example~\ref{ex:impact_uncertainty} and
\Cref{fig:gbm_histogram_a_cec_hjb} it appears that the relative
improvement in profits by pricing according to $\alpha^B(t)$, rather than the heuristic
policy $\alpha^D(t)$, is negligible.
Further numerical experiments for other values of $\sigma$ strengthen
these results. See \Cref{tbl:profit_hjb_cec_statistics} for summary
statistics of the relative difference in profits
$1-P(\alpha^D)/P(\alpha^B)$.
The relative improvement of following the strategy $\alpha^B(t)$
increases with $\sigma$, however the standard deviations are all on
the order of \SIrange{0.01}{0.1}{\percent}.
\begin{table}[htb]
  \centering
  \begin{tabular}{lrrrrr}
    \toprule
    $\sigma$ & mean & std & $Q_{0.05}$ & $Q_{0.5}$ & $Q_{0.95}$\\
    \midrule
    0.05& 0.0 & 0.4 &-0.5&-0.2&0.9\\
    0.1&0.2& 0.6&-0.4&-0.1&1.4\\
    0.2&0.3& 0.7&-0.3&0.0&1.8\\
    0.4&0.4& 1.2&-0.3&0.1&1.8\\
    \midrule
             & $\times 10^{-3}$&$\times 10^{-3}$&$\times 10^{-3}$&$\times 10^{-3}$&$\times 10^{-3}$\\
    \bottomrule
  \end{tabular}
  \caption{Summary statistics of the relative profit difference
    $1-P(\alpha^D)/P(\alpha^B)$ from the model in
    Example~\ref{ex:impact_uncertainty} with different levels of
    uncertainty $\sigma$.
    The headings $Q_z$ denote the $z$-quantile of the distribution.
  }\label{tbl:profit_hjb_cec_statistics}
\end{table}

\section{Extensions of the pricing problem}\label{sec:extensions}
For our one-product system with multiplicative parameter dynamics with uncertainty,
two natural extensions to the pricing problem are: (i) To incorporate
other forms of uncertainty, and (ii) formulate the problem for risk-averse
decision makers. For example, the costs of handling unsold stock may
not be known a priori, and risk aversion may be modelled from an
expected utility viewpoint.
These extensions increase the input dimensions of the corresponding
HJB equations, which we present without further discussion. We again
refer readers interested in the
derivation of the HJB equations to~\cite{pham2009continuous}.

\subsection{Other forms of uncertainty}
In the prior sections, the source of randomness in the system has come from
the multiplicative term $G(t)$, modelled as a geometric Brownian
motion martingale. A different non-negative stochastic process may be more appropriate, and the
choice of dynamics can be guided by existing sales data.
More generally, the demand function $q(a)$ depends on multiple
parameters that exhibit different levels of
uncertainty and dynamics in time. Another parameter in the pricing
problem is the unit cost $C$, where its value at the terminal time may
depend on factors unknown at times $t<1$.
Let $\theta(t)\in\mathbb R^n$ denote the vector of parameters that are
relevant to the problem, and say the demand function  and the unit cost
depend explicitly on $\theta$. We write $q(a;\theta)$ and $C(\theta)$
for this dependence.
Within the diffusion based stochastic framework, we can model the
dynamics of $\theta(t)$ with functions $b(t,\theta)\in\mathbb R^n$,
$\sigma(t,\theta)\in\mathbb R^{n\times p}$, and a vector-valued, uncorrelated,
Brownian motion $W(t)\in\mathbb R^p$.
We assume that $\theta(t)$ does not depend on the pricing policy
$\alpha(t)$ or the remaining
stock $S^\alpha(t)$. Thus,
the system for the pricing problem is described by the SDE,
\begin{align}
  d\theta(t)&=b(t,\theta(t))\,\mathrm{d}t + \sigma(t,\theta(t))\,\mathrm{d}W(t),\\
  \mathrm{d}S^\alpha(t)&=-q(\alpha(t);\theta(t)).
\end{align}

Let $\Theta\in\mathbb R^n$ denote the state space for
$\theta(t)$.
Let $D_\theta v$ and $D_\theta^2v$ denote the gradient and Hessian of
$v(t,s,\theta)$ with respect to $\theta$. We denote the transpose
operator with a superscript asterisk, and introduce the
volatility matrix $\Sigma(t,\theta) = \sigma(t,\theta){\sigma(t,\theta)}^*$.
The HJB approach for the pricing problem is then to find
$v:{[0,1]}^2\times\Theta\to\mathbb R$ which satisfies
\begin{align}
  v_t+{b(t,\theta)}^* D_\theta v
  + {\textstyle\frac{1}{2}}\mbox{tr}\left( \Sigma(t,\theta)
  D_\theta^2v \right)
  +\max_{a}\{q(a;\theta)(a-v_s)\}&=0,\\
  v(1,s,\theta) &= -C(\theta)s.
\end{align}
Additional boundary conditions may be necessary, depending on
$\theta(t)$. For example, the boundary condition for $g=0$ in our
multiplicative model from the previous sections.
The value function and the pricing function now depend
explicitly on each element in $\theta$, thus increasing the dimension
of the corresponding HJB equation. Efficient algorithms for solving
high-dimensional PDEs of this form exist. See, for example, the
multigrid preconditioning approach by~\cite{reisinger2017boundary}, or
the splitting into a sequence of lower-dimensional PDEs by~\cite{reisinger2017finite}.
In higher dimensions, it is of even greater importance to critically
balance computational cost with the suboptimality of
approximations. As \Cref{sec:stochastic_hjb} shows, approximate
policies can perform well. This underscores the importance of assessing whether
new sources of randomness in the model have significant impact on the
objective and the optimal policy.

\subsection{Expected utility risk aversion}
One may argue that for a retailer as a whole, an assumption of a
risk-neutral decision maker is valid. For individual product managers,
whose performance is evaluated over shorter time horizons, a degree of
risk aversion can be preferable from their point of view. Formulations and investigations of the
impact of risk aversion on pricing policies is also noted as an
interesting area of research by~\cite{bitran2003overview}.
For investigations of the expected utility problem with Poisson process
demand, see, for example, \citet{lim2007relative,feng2008risk}.

In this section, we assume that the decision maker is evaluated based on
the performance of the total profit from selling a product over the
time interval $[0,1]$. Then, the pricing decision at time $t$ may also depend
on how much revenue has been accrued at that time. So we
introduce a state variable $R^\alpha(t)$ representing the accrued
revenue at time $t$, with dynamics
$dR^\alpha(t)=\alpha(t)q(\alpha(t);\theta(t))\,\mathrm{d}t$.
We consider risk aversion based on
an expected utility-maximizing decision maker.
Given a utility function $U(x)$, the pricing problem is then to
find a pricing policy $\alpha$ that maximizes the expected utility
$\mathbb E[U(R^\alpha(1)-C(\theta(1))S^\alpha(1))]$.
For this utility function, the value function is defined as
\begin{equation}
  v(t,s,\theta,r)=
  \max_{\alpha\in\mathcal A}
  \mathbb E_t\left[
    U\left(R^\alpha(1)-C(\theta(1))S^\alpha(1)\right)
  \right].
\end{equation}
The corresponding HJB problem is then to find
$v:{[0,1]}^2\times \Theta\times{[0,\infty)}\to\mathbb R$ that solves
\begin{align}\label{eq:hjb_utility}
  v_t+{b(t,\theta)}^* D_\theta v
  + {\textstyle\frac{1}{2}}\mbox{tr}\left( \Sigma(t,\theta)
  D_\theta^2v \right)
  +\max_{a}\{q(a;\theta)(av_r-v_s)\}&=0,\\
  v(1,s,\theta,r) &= U(r-C(\theta)s),
\end{align}
plus additional boundary conditions. The impact of the risk aversion
on pricing decisions appears in the maximization term in the HJB
equation, $\max_a\{q(a;\theta)(av_r-v_s)\}$. The $v_r$ term represents
the relative importance of accruing more revenue to the utility of
selling more stock, as represented by the $v_s$ term.

\section{Conclusion}\label{sec:conclusion}
This article focuses on continuum approximations for dynamic
pricing problems under uncertainty.
Most of the existing literature on continuous time dynamic pricing for
revenue management is based on Poisson
processes. This is suitable for many applications,
however for large retailers, approximating the number of sales and
stock as a continuum can simplify calculation of pricing rules.
We present an approach for modelling the sales as a continuous time
dynamical system, where the uncertainty in demand arises from
stochastic processes.
An advantage over the Poisson process model is that this approach
allows us to directly model the demand volatility.
Under this model, we consider a pricing problem where the retailer
aims to deplete  inventory of a product at maximum profit.
By formulating the problem as a stochastic optimal control problem,
we can express the optimal pricing policy in terms of the solution to
a nonlinear PDE.\@
For linear and exponential demand functions, we find closed-form
expressions for the pricing policy when the system is deterministic.
It turns out that, for a risk-neutral decision maker, the
deterministic pricing policy is a near-optimal heuristic for systems
with demand uncertainty.
Numerical errors in calculating the optimal pricing policy
may in fact result in lower profits on average than with the heuristic
pricing policy.

There are two topics of particular interest for future study.
The first is to understand why demand uncertainty has such a
small effect on the optimal pricing policy for risk-neutral decision
makers, and whether constraints such as requiring monotone-in-time pricing
policies may increase this effect. Second, a case study of the continuum
model framework for multiple products and time-dependent demand is
needed, in
order to understand how well this approach can scale to revenue
management implementations for retailers.

\ifanon
\else
\section{Funding}
This work was supported by the EPSRC
Centre For Doctoral Training in Industrially Focused Mathematical
Modelling (EP/L015803/1) in collaboration with dunnhumby Limited.

\section{Acknowledgements}
The author would like to thank his supervisors Jeff Dewynne and
Chris Farmer for their input and discussions leading to this
work. The input from Jon Chapman and Jeff Dewynne on the asymptotic
analysis of the linear-demand model has been invaluable.
I also thank two anonymous reviewers and the editor for their useful
suggestions and criticism.
\fi

\appendix
\section*{Appendix: Asymptotic analysis}
We carry out the asymptotic analysis of the pricing function
that arises from the linear demand HJB
equation~\eqref{eq:hjb_impact_uncertainty} under the assumptions of
\Cref{sec:stochastic_hjb}.
Further, we assume that the value and pricing functions are sufficiently
differentiable to carry out the operations below. The smoothing effect
of the diffusion in the PDEs when $\sigma>0$ justifies this assumption.
Result~\ref{res:linear_asymptotics} in \Cref{sec:stochastic_hjb}
discusses the interpretation of the following results.

\begin{figure}[hbt]
  \centering
  \includegraphics[width=0.5\textwidth,axisratio=1.5]{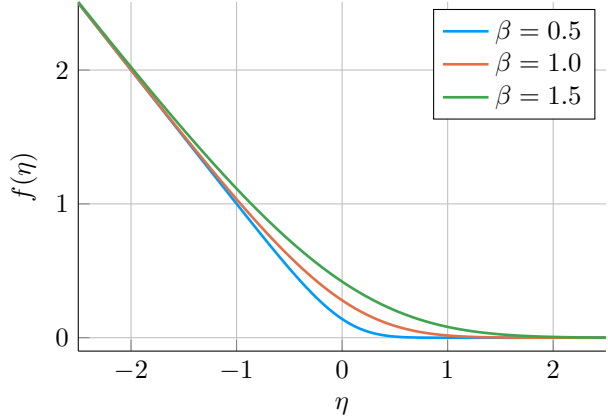}
  \caption{Numerically computed solutions of~\eqref{eq:bvp_ode_leading}.
    Note the similarity to the optimal pricing function curves zoomed
    in near the kink in \Cref{fig:hjb_cec_zoom}.
  }\label{fig:inner_leading_f}
\end{figure}

Let us first reduce the dimension of the problem with a similarity
transform working in reverse time, and then present a PDE satisfied by the pricing function.
Consider the transformations $\xi(s,g) = s/g$, $\tau(t)=1-t$, and
$v(t,s,g)=g\phi(\tau(t),\xi(s,g))$, which satisfies
\begin{align}\label{eq:hjb_value_transform}
  \phi_\tau &= \frac{1}{2}\sigma^2\xi^2\phi_{\xi\xi} +
           \frac{1}{4}{(q_1 - \phi_\xi)}^2 \\
  \phi(0,\xi) &= -C\xi\\
  \phi(\tau,0)& = 0,\\
  \phi_\xi(\tau,\infty) &= -C.
\end{align}
In the new coordinates, define the pricing function in terms of the
function $\psi(\tau,\xi)$ so that $a^B(t,s,g)=\psi(\tau(t),\xi(s,g))$.
Then $\psi = (q_1+\phi_\xi)/2$ satisfies the following PDE
that arises from differentiating~\eqref{eq:hjb_value_transform} with
respect to $\xi$.
\begin{align}\label{eq:hjb_policy_transform}
  \psi_\tau&= \frac{1}{2}\sigma^2 \xi^2
          \psi_{\xi\xi} +
          (\sigma^2\xi - q_1 + \psi)\psi_\xi,\\
  \psi(\tau,0) &= q_1,\\
  \psi(\tau,\infty) &= \delta:=\max(0,(q_1-C)/2).
\end{align}
When $\sigma=0$, we know from \Cref{sec:deterministic_hjb} that
the viscosity solution $\psi^{(0)}$ is
\begin{equation}
  \psi^{(0)}(\tau,\xi)=\begin{cases}
    q_1-\xi/\tau, &\text{if } 0\leq \xi\leq
    \beta \tau,\\
    \delta,&\text{if } \hfill \xi > \beta \tau,
  \end{cases}
  \qquad \textnormal{where } \beta:=\min(q_1,(q_1+C)/2).
\end{equation}
This leads us to consider an inner layer asymptotic analysis near
the kink $\xi = \beta \tau$. Zoom in near the kink using the coordinates
$x(\tau,\xi) = (\xi-\beta \tau)/\sigma$ and write $\psi(\tau,\xi) = \sigma
u(\tau,x(\tau,\xi)) + \delta$. In the new
coordinates~\eqref{eq:hjb_policy_transform} becomes
\begin{align}
  u_\tau 
  &=\frac{1}{2}{(\sigma x + \beta \tau)}^2u_{xx}
    +\sigma (\sigma x+\beta \tau)u_x  
    + uu_x,&\tau>0,\,x&> -\beta \tau/\sigma,
\end{align}
with matching conditions $\lim_{x \to -\infty} u_x(\tau,x) = -1/\tau$ and
$\lim_{x\to \infty}u(\tau,x) = 0$.
The leading order equation for $u = u^{(0)}+\sigma u^{(1)}+\sigma^2 u^{(2)}+\cdots$ is therefore
\begin{equation}
  u^{(0)}_\tau = \frac{1}{2}{(\beta \tau)}^2u^{(0)}_{xx} + u^{(0)}u^{(0)}_x.
\end{equation}
This equation has a similarity solution using the transformations
$\eta(\tau,x)= x/\tau^{3/2}$ and $u^{(0)}(\tau,x) = \tau^{1/2}f(\eta(\tau,x))$.
The function $f(\eta)$ must satisfy the boundary-value ODE
\begin{align}\label{eq:bvp_ode_leading}
  \beta^2f^{\prime\prime} + (3\eta + 2f)f^\prime - f
  &= 0,\\
  f&\to 0 &\textnormal{as } \eta&\to \infty,\\
  f^\prime &\to -1&\textnormal{as } \eta&\to -\infty.
\end{align}
\Cref{fig:inner_leading_f} shows numerically computed solutions to the
ODE for different values of $\beta$.

We briefly note that an asymptotic expansion of
$\phi=\phi^{(0)}+\sigma\phi^{(1)}+\sigma^2\phi^{(2)}+\cdots$
in the outer region $\xi<\beta t$ reveals that $\phi^{(1)}=0$ and
$\phi^{(2)}(t,\xi)=-\xi^2/2$. It follows that
$\psi(t,\xi)\approx \psi^{(0)}(t,\xi)-\sigma^2 \xi/2$ in this region.

\bibliographystyle{agsm}
\bibliography{references}

\end{document}